\documentclass[12pt]{article}
\usepackage{amssymb}
\usepackage{amsfonts}
\usepackage{amsmath}
\usepackage{amsbsy}
\newcommand{\be}{\begin{equation}}
\newcommand{\ee}{\end{equation}}
\newcommand{\bea}{\begin{eqnarray}}
\newcommand{\eea}{\end{eqnarray}}
\newcommand{\ba}{\begin{array}}
\newcommand{\ea}{\end{array}}

\newcommand{\bc}{\begin{center}}
\newcommand{\ec}{\end{center}}
\newcommand{\ben}{\begin{enumerate}}
\newcommand{\een}{\end{enumerate}}
\newcommand{\bfi}{\begin{figure}}
\newcommand{\efi}{\end{figure}}

\newcommand{\bq}{\begin{quote}}
\newcommand{\eq}{\end{quote}}
\newcommand{\bqu}{\begin{quotation}}
\newcommand{\equ}{\end{quotation}}
\newenvironment{emphit}{\begin{itemize}}{\end{itemize}}
\newcommand{\bemp}{\begin{emphit}}
\newcommand{\eemp}{\end{emphit}}

\newcommand{\bt}{\begin{tabular}}
\newcommand{\et}{\end{tabular}}

\newtheorem{myth}{Theorem}[section]
\newtheorem{mylem}{Lemma}[section]
\newtheorem{mycor}{Corollary}[section]

\newtheorem{mydef}{Definition}[section]
\newtheorem{myrem}{Remark}[section]

\begin{document}
\date{}
\title{On A Cryptographic Identity In Osborn Loops \footnote{2000
Mathematics Subject Classification. Primary 20NO5 ; Secondary 08A05}
\thanks{{\bf Keywords and Phrases :} Osborn loops, cryptography}}
\author{T. G. Jaiy\'e\d ol\'a\thanks{All correspondence to be addressed to this author.} \\
Department of Mathematics,\\
Obafemi Awolowo University,\\
Ile Ife 220005, Nigeria.\\
jaiyeolatemitope@yahoo.com\\tjayeola@oauife.edu.ng \and
J. O. Ad\'en\'iran \\
Department of Mathematics,\\
University of Abeokuta, \\
Abeokuta 110101, Nigeria.\\
ekenedilichineke@yahoo.com\\
adeniranoj@unaab.edu.ng} \maketitle
\begin{abstract}
This study digs out some new algebraic properties of an Osborn loop
that will help in the future to unveil the mystery behind the middle
inner mappings $T_{(x)}$ of an Osborn loop. These new algebraic
properties, will open our eyes more to the study of Osborn loops
like CC-loops which has received a tremendious attention in this
$21^\textrm{st}$ and VD-loops whose study is yet to be explored. In
this study, some algebraic properties of non-WIP Osborn loops have
been investigated in a broad manner. Huthnance was able to deduce
some algebraic properties of Osborn loops with the WIP i.e universal
weak WIPLs. So this work exempts the WIP. Two new loop identities,
namely left self inverse property loop(LSIPL) identity and right
self inverse property loop(RSLPL) are introduced for the first time
and it is shown that in an Osborn loop, they are equivalent. A
CC-loop is shown to be power associative if and only if it is a
RSLPL or LSIPL. Among the few identities that have been established
for Osborn loops, one of them is recognized and recommended for
cryptography in a similar spirit in which the cross inverse property
has been used by Keedwell following the fact that it was observed
that Osborn loops that do not have the LSIP or RSIP or 3-PAPL or
weaker forms of inverse property, power associativity and
diassociativity to mention a few, will have cycles(even long ones).
These identity is called an Osborn cryptographic identity(or just a
cryptographic identity).
\end{abstract}

\section{Introduction}
\paragraph{}
Let $L$ be a non-empty set. Define a binary operation ($\cdot $) on
$L$ : If $x\cdot y\in L$ for all $x, y\in L$, $(L, \cdot )$ is
called a groupoid. If the system of equations ;
\begin{displaymath}
a\cdot x=b\qquad\textrm{and}\qquad y\cdot a=b
\end{displaymath}
have unique solutions for $x$ and $y$ respectively, then $(L, \cdot
)$ is called a quasigroup. Furthermore, if there exists a unique
element $e\in L$ called the identity element such that for all $x\in
L$, $x\cdot e=e\cdot x=x$, $(L, \cdot )$ is called a loop. We write
$xy$ instead of $x\cdot y$, and stipulate that $\cdot$ has lower
priority than juxtaposition among factors to be multiplied. For
instance, $x\cdot yz$ stands for x(yz). For each $x\in L$, the
elements $x^\rho =xJ_\rho ,x^\lambda =xJ_\lambda\in L$ such that
$xx^\rho=e=x^\lambda x$ are called the right, left inverses of $x$
respectively. $x^{\lambda^i}=(x^\lambda )^\lambda$ and
$x^{\rho^i}=(x^\rho )^\rho$ for $i\ge 1$. $L$ is called a weak
inverse property loop (WIPL) if and only if it obeys the weak
inverse property (WIP);
\begin{displaymath}
xy\cdot z=e~\textrm{implies}~x\cdot yz=e~\textrm{for all}~x,y,z\in L
\end{displaymath}
while $L$ is called a cross inverse property loop (CIPL) if and only
if it obeys the cross inverse property (CIP);
\begin{displaymath}
xy\cdot x^\rho=y.
\end{displaymath}
The triple $\alpha =(A,B,C)$ of bijections on a loop $(L,\cdot )$ is
called an autotopism of the loop if and only if
\begin{displaymath}
xA\cdot yB=(x\cdot y)C~\textrm{for all}~x,y\in L.
\end{displaymath}
Such triples form a group $AUT(L,\cdot )$ called the autotopism
group of $(L,\cdot )$. In case the three bijections are the same i.e
$A=B=C$, then any of them is called an automorphism and the group
$AUM(L,\cdot )$ which such forms is called the automorphism group of
$(L,\cdot )$. For an overview of the theory of loops, readers may
check \cite{phd3,phd41,phd39,phd49,phd42,phd75}.

Osborn \cite{phd89}, while investigating the universality of WIPLs
discovered that a universal WIPL $(G,\cdot )$ obeys the identity
\begin{equation}\label{eq:1}
yx\cdot (zE_y\cdot y)=(y\cdot xz)\cdot y~\textrm{for all}~x,y,z\in G
\end{equation}
\begin{displaymath}
\textrm{where}~\theta_y=L_yL_{y^\lambda}=R_{y^\rho}^{-1}R_y^{-1}=L_yR_yL_y^{-1}R_y^{-1}.
\end{displaymath}
A loop that necessarily and sufficiently satisfies this identity is
called an Osborn loop.

\paragraph{}
Eight years after Osborn's \cite{phd89} 1960 work on WIPL, in 1968,
Huthnance  Jr. \cite{phd44} studied the theory of generalized
Moufang loops. He named a loop that obeys (\ref{eq:1}) a generalized
Moufang loop and later on in the same thesis, he called them
M-loops. On the other hand, he called a universal WIPL an Osborn
loop and this same definition was adopted by Chiboka \cite{phd96}.
Basarab \cite{phd148,phd46,phd137} and Basarab and Belioglo
\cite{phd170} dubbed a loop $(G,\cdot )$ satisfying any of the
following equivalent identities an Osborn loop:
\begin{equation}\label{eq:3}
OS_2~:~x(yz\cdot x)=(x^\lambda\backslash y)\cdot  zx
\end{equation}
\begin{equation}\label{eq:4}
OS_3~:~(x\cdot yz)x=xy\cdot (zE_x^{-1}\cdot x)
\end{equation}
\begin{displaymath}
\textrm{where}~E_x=R_xR_{x^\rho}=(L_xL_{x^\lambda})^{-1}=R_xL_xR_x^{-1}L_x^{-1}~\textrm{for
all}~x,y,z\in G
\end{displaymath} and the binary operations '$\backslash$' and
'$/$' are respectively defines as ; $z=x\cdot y$ if and only if
$x\backslash z=y$ if and only if $z/y=x$ for all $x,y,z\in G$.

It will look confusing if both Basarab's and Huthnance's definitions
of an Osborn loop are both adopted because an Osborn loop of Basarab
is not necessarily a universal WIPL(Osborn loop of Huthnance). So in
this work, Huthnance's definition of an Osborn loop will be dropped
while we shall stick to that of Basarab which was actually adopted
by M. K. Kinyon \cite{phd33} who revived the study of Osborn loops
in 2005 at a conference tagged "Milehigh Conference on Loops,
Quasigroups and Non-associative Systems" held at the University of
Denver, where he presented a talk titled "A Survey of Osborn Loops".

Let $t=x^\lambda\backslash y$ in $OS_2$, then $y=x^\lambda t$ so
that we now have an equivalent identity
$$x[(x^\lambda y)z\cdot x]=y\cdot  zx.$$
Huthnance \cite{phd44} was able to deduce some properties of $E_x$
relative to (\ref{eq:1}). $E_x=E_{x^\lambda}=E_{x^\rho}$. So, since
$E_x=R_xR_{x^\rho}$, then $E_x=E_{x^\lambda}=R_{x^\lambda}R_{x}$ and
$E_x=(L_{x^\rho}L_x)^{-1}$. So, we now have the following equivalent
identities defining an Osborn loop.
\begin{equation}\label{eq:3.1}
OS_2~:~x[(x^\lambda y)z\cdot x]=y\cdot  zx
\end{equation}
\begin{equation}\label{eq:4.1}
OS_3~:~(x\cdot yz)x=xy\cdot [(x^\lambda \cdot xz)\cdot x]
\end{equation}

\begin{mydef}
A loop $(Q,\cdot )$ is called:
\begin{description}
\item[(a)] a 3 power associative property loop(3-PAPL) if and only if $xx\cdot
x=x\cdot xx$ for all $x\in Q$.
\item[(b)] a left self inverse property loop(LSIPL) if and only if
$x^\lambda\cdot xx=x$ for all $x\in Q$.
\item[(c)] a right self inverse property loop(RSIPL) if and only if
$xx\cdot x^\rho =x$ for all $x\in Q$.
\end{description}
\end{mydef}

The identities describing the most popularly known varieties of
Osborn loops are given below.
\begin{mydef}
A loop $(Q,\cdot )$ is called:
\begin{description}
\item[(a)] a VD-loop if and only if
\begin{displaymath}
(\cdot )_x=(\cdot )^{L_x^{-1}R_x}\qquad\textrm{and}\qquad {}_x(\cdot
)=(\cdot )^{R_x^{-1}L_x} \end{displaymath} i.e $R_x^{-1}L_x\in
PS_\lambda (Q,\cdot )$ with companion $c=x$ and $L_x^{-1}R_x\in
PS_\rho (Q,\cdot )$ with companion $c=x$ for all $x\in Q$ where
$PS_\lambda (Q,\cdot )$ and $PS_\rho (Q,\cdot )$ are respectively
the left and right pseudo-automorphism groups of $Q$.\qquad Basarab
\cite{phd137}
\item[(b)] a Moufang loop if and only if the identity
\begin{displaymath}
(xy)\cdot (zx)=(x\cdot yz)x
\end{displaymath}
holds in $Q$.
\item[(c)] a conjugacy closed loop(CC-loop) if and only if the
identities
\begin{displaymath}
x\cdot yz=(xy)/x\cdot xz\qquad\textrm{and}\qquad zy\cdot x=zx\cdot
x\backslash(yx)
\end{displaymath}
hold in $Q$.
\item[(d)] a universal WIPL if and only if the identity
\begin{displaymath}
x(yx)^\rho=y^\rho\qquad\textrm{or}\qquad(xy)^\lambda x=y^\lambda
\end{displaymath}
holds in $Q$ and all its isotopes.
\end{description}
\end{mydef}
All these three varieties of Osborn loops and universal WIPLs are
universal Osborn loops. CC-loops and VD-loops are G-loops. G-loops
are loops that are isomorphic to all their loop isotopes. Kunen
\cite{phd185} has studied them.

In the multiplication group ${\cal M}(G,\cdot )$ of a loop $(G,\cdot
)$ are found three important permutations, namely, the right, left
and middle inner mappings $R_{(x,y)}=R_xR_yR_{xy}^{-1}$,
$L_{(x,y)}=L_xL_yL_{yx}^{-1}$ and $T_{(x)}=R_xL_x^{-1}$ respectively
which form the right inner mapping group $\textrm{Inn}_\lambda(G)$,
left inner mapping group $\textrm{Inn}_\rho (G)$ and the middle
inner mapping $\textrm{Inn}_\mu (G)$. In a Moufang loop $G$,
$R_{(x,y)},L_{(x,y)},T_{(x)}\in PS_\rho(G)$ with companions
$(x,y),(x^{-1},y^{-1}),x^{-3}\in G$ respectively.
\begin{myth}\label{1:3.02}(Kinyon \cite{phd33})

Let $G$ be an Osborn loop. $R_{(x,y)}\in PS_\rho(G)$ with companion
$(xy)^\lambda (y^\lambda\backslash x)$ and $L_{(x,y)}\in
PS_\lambda(G)~\forall~x,y\in G$. Furthermore,
$R_{(x,y)}^{-1}=[L_{y^\rho}^{-1},R_x^{-1}]=L_{(y^\lambda
,x^\lambda)}~\forall~x,y\in G$.
\end{myth}
The second part of Theorem~\ref{1:3.02} is trivial for Moufang
loops. For CC-loops, it was first observed by Dr\'apal and then
later by Kinyon and Kunen \cite{phd47}.
\begin{myth}\label{1:3.05}
Let $G$ be an Osborn loop.
$\textrm{Inn}_\rho(G)=\textrm{Inn}_\lambda(G)$.
\end{myth}
Still mysterious are the middle inner mappings $T_{(x)}$ of an
Osborn loop. In a  Moufang loop, $T_{(x)}\in PS_\rho$ with a
companion $x^{-3}$ while in a CC-loop, $T_{(x)}\in PS_\lambda$ with
companion $x$. Kinyon \cite{phd33}, possess a question asking of
what can be said in case of an arbitrary Osborn loop.
\begin{myth}\label{1:3.050}(Kinyon \cite{phd33})

In an Osborn loop $G$ with centrum $C(G)$ and center $Z(G)$:
\begin{enumerate}
\item If $T_{(a)}\in AUM(G)$, then $a\cdot aa=aa\cdot a\in N(G)$. Thus, for all $a\in C(G)$, $a^3\in Z(G)$.
\item If $(xx)^\rho =x^\rho x^\rho$ holds, then
$x^{\rho\rho\rho\rho\rho\rho}=x$ for all $x\in G$.
\end{enumerate}
\end{myth}
Some basic loop properties such as flexibility, left alternative
property(LAP), left inverse property(LIP), right alternative
property(RAP), right inverse property(RIP), anti-automorphic inverse
property(AAIP) and the cross inverse property(CIP) have been found
to force an Osborn loop to be a Moufang loop. This makes the study
of Osborn loops more challenging and care must be taking not to
assume any of these properties at any point in time except the WIP,
automorphic inverse property and some other generalizations of the
earlier mentioned loop properties(LAP, LIP, e.t.c.).

\begin{mylem}\label{1:4.1}
An Osborn loop that is flexible or which has the LAP or RAP or LIP
or RIP or AAIP is a Moufang loop. But an Osborn loop that is
commutative or which has the CIP is a commutative Moufang loop.
\end{mylem}
\begin{myth}\label{0:2}(Basarab, \cite{phd46})

If an Osborn loop is of exponent 2, then it is an abelian group.
\end{myth}

\begin{myth}\label{1:3.13}(Huthnance \cite{phd44})

Let $G$ be a WIPL. $G$ is a universal WIPL if and only if $G$ is an
Osborn loop.
\end{myth}

\begin{mylem}\label{3:13}(Lemma~2.10, Huthnance \cite{phd44})

Let $L$ be a WIP Osborn loop. If $a=x^\rho x$, then for all $x\in
L$:
\begin{displaymath}
xa=x^{\lambda^2},~ax^\lambda =x^\rho,~x^\rho
a=x^\lambda,~ax=x^{\rho^2},~xa^{-1}=ax,~a^{-1}x^\lambda =x^\lambda
a,~a^{-1}x^\rho =x^\rho a.
\end{displaymath}
or equivalently
\begin{displaymath}
J_\lambda~:~x\mapsto~x\cdot x^\rho x,~J_\rho~:~x\mapsto~x^\rho
x\cdot x^\lambda,~J_\lambda~:~x\mapsto~x^\rho\cdot x^\rho x
,~J_\rho^2~:~x\mapsto~x^\rho x\cdot x,
\end{displaymath}
\begin{displaymath}
x(x^\rho x)^{-1}=(x^\rho x)x,~(x^\rho x)^{-1}x^\lambda =x^\lambda
()x^\rho x,~(x^\rho x)^{-1}x^\rho =x^\rho (x^\rho x).
\end{displaymath}
\end{mylem}
The aim of this study is to dig out some new algebraic properties of
an Osborn loop that will help in the future to unveil the mystery
behind the middle inner mappings $T_{(x)}$ of an Osborn loop. These
new algebraic properties, will open our eyes more to the study of
Osborn loops like CC-loops, introduced by Goodaire and Robinson
\cite{phd91,phd48}, whose algebraic structures have been studied by
Kunen \cite{phd78} and some recent works of Kinyon and Kunen
\cite{phd36,phd47}, Phillips et. al. \cite{phd35}, Dr\'apal
\cite{phd37,phd38,phd98,phd107}, Cs\"org\H o et. al.
\cite{phd106,phd104,phd108} and VD-loops whose study is yet to be
explored. In this study, the algebraic properties of non-WIP Osborn
loops have been investigated in a broad manner. Huthnance
\cite{phd44} was able to deduce some algebraic properties of Osborn
loops with the WIP i.e universal WIPLs. So this work exempts the
WIP. Two new loop identities, namely left self inverse property
loop(LSIPL) identity and right self inverse property loop(RSLPL) are
introduced for the first time and it is shown that in an Osborn
loop, they are equivalent. A CC-loop is shown to be power
associative if and only if it is a RSLPL or LSIPL. Among the few
identities that have been established for Osborn loops, one of them
is recognized and recommended for cryptography in a similar spirit
in which the cross inverse property has been used by Keedwell
following the fact that it was observed that Osborn loops that do
not have the LSIP or RSIP or 3-PAPL or weaker forms of inverse
property, power associativity and diassociativity to mention a few,
will have cycles(even long ones). These identity is called an Osborn
cryptographic identity(or just a cryptographic identity).

\section{Main Results}
\subsection{Some Algebraic Properties Of Osborn Loops}
\begin{myth}\label{3:16}
Let $(L,\cdot )$ be a loop. $L$ is an Osborn loop if and only if
$(L_{x^\lambda}, R_x^{-1},L_x^{-1}R_x^{-1})\in AUT(L)$. Hence for
all $x,y,z\in L$:
\begin{enumerate}
\item $(L_{(x)},L_{(x)},L_xT_{(x)}R_x^{-1})\in AUT(L)$ for some
$L_{(x)}\in \textrm{Inn}_\lambda(L)$.
\item
\begin{enumerate}
\item $T_{(x)}~:~x\mapsto~[(x^\lambda\cdot xy)(x^\lambda\cdot xy^\rho
)]x$.
\item $T_{(x)}~:~x\mapsto~[(x^\lambda\cdot xz^\lambda )(x^\lambda\cdot xz)]x$.
\item $T_{(x)}~:~y\mapsto~x^\lambda y\cdot x$ i.e
$T_{(x)}:=L_{x^\lambda }R_x$.
\end{enumerate}
\item $yx=x(x^\lambda y\cdot x)$ i.e $R_x=L_{x^\lambda}R_xL_x$.
\item $(x^\lambda\cdot xy)(x^\lambda\cdot xy^\rho )=(x^\lambda\cdot xz^\lambda)(x^\lambda xz
)=e$.
\end{enumerate}
\end{myth}
{\bf Proof}\\
By $OS_2$, $L$ is an Osborn loop if and only if
$(L_{x^\lambda},R_x^{-1},L_x^{-1}R_x^{-1})\in AUT(L)$. By $OS_3$,
$L$ is an Osborn loop if and only if
$(L_x,L_xL_{x^\lambda}R_x,L_xR_x)\in AUT(L)$.
\begin{enumerate}
\item Hence, $(L_{(x)},L_{(x)},L_xT_{(x)}R_x^{-1})\in AUT(L)$ where $L_{(x)}=L_xL_{x^\lambda}\in \textrm{Inn}_\lambda(L)$.

The autotopism $(L_{(x)},L_{(x)},L_xT_{(x)}R_x^{-1})$ implies
$[(x^\lambda\cdot xy)\cdot (x^\lambda\cdot xz)]x=(x\cdot
yz)T_{(x)}$.
\item
\begin{enumerate}
\item So with $z=y^\rho$, $xT_{(x)}=[(x^\lambda\cdot xy)\cdot (x^\lambda\cdot xy^\rho)]x$.
\item Similarly, with $y=z^\lambda$, $xT_{(x)}=[(x^\lambda\cdot xz^\lambda )\cdot (x^\lambda\cdot xz)]x$.
\item With $y=e$ or $z=e$, $(xz)T_{(x)}=(x^\lambda\cdot xz)x$ which
implies that $yT_{(x)}=(x^\lambda\cdot y)x$ or
$zT_{(x)}=(x^\lambda\cdot z)x$ respectively.
\end{enumerate}
\item Recall that $T_{(x)}=R_xL_x^{-1}$. Using this and
$T_{(x)}=L_{x^\lambda}R_x$, $R_x=L_{x^\lambda}R_xL_x$.
\item Observe that $xT_{(x)}=x$, so by (b)i. and (b)ii. the claim is
true.
\end{enumerate}

\begin{mylem}\label{3:17}
Let $(L,\cdot )$ be an Osborn loop. The following are true.
\begin{enumerate}
\item $(x^\lambda\cdot xy)^\rho =x^\lambda\cdot xy^\rho$,~$(x^\lambda\cdot xy^\rho )^\lambda =(x^\lambda\cdot xy^\lambda
)^\rho$.
\item $J_\rho~:~x~\mapsto x^\lambda x^\lambda\cdot x$,~
$J_\rho^2~:~x\mapsto xx\cdot x^\rho$,~
$J_\lambda~:~x\mapsto~(x^\lambda )^\lambda\cdot x^\lambda
x^\rho$,~$J_\lambda^2~:~x\mapsto~x^\lambda\cdot
xx$,~$J_\lambda~:~x\mapsto~(x^\lambda\cdot
xx^\lambda)^2(x^\lambda\cdot xx)$,~$J_\lambda~:~x\mapsto~(x^\lambda
x^\lambda\cdot x)^\lambda (x^\lambda x^\lambda\cdot
x)^2$,~$J_\lambda^3~:~x~\mapsto x^\lambda\cdot xx^\lambda$
\item \begin{displaymath}
x^\lambda\cdot xx^{\rho^2}=x,\qquad (x\cdot x^\rho x^\rho )^\lambda
=x\cdot x^\rho x=(x\cdot x^\rho x^\lambda )^\rho,
\end{displaymath}
\begin{displaymath}
(x^\lambda\cdot xx)^\lambda =x^\lambda\cdot xx^\lambda,\qquad
x^{\lambda^3}\cdot x^{\lambda^2}x=x^\lambda\cdot xx,
\end{displaymath}
\begin{displaymath}
(x^{\lambda^2}\cdot x^\lambda x^\rho )^{\lambda^2}\cdot
(x^{\lambda^2}\cdot x^\lambda x^\rho )^\lambda (x^{\lambda^2}\cdot
x^\lambda x^\rho )^\rho=x^\lambda\cdot xx,\qquad (x^\lambda\cdot
xx^\lambda )^2(x^\lambda\cdot xx)=x^{\lambda^2}\cdot x^\lambda
x^\rho,
\end{displaymath}
\begin{displaymath}
(x^\lambda x^\lambda\cdot x)^\lambda (x^\lambda x^\lambda\cdot x)^2
=(x^\lambda\cdot xx^\lambda )^2(x^\lambda\cdot xx),
\end{displaymath}
\begin{displaymath}
(x^{\lambda^2}\cdot x^\lambda x^\rho )^\lambda (x^{\lambda^2}\cdot
x^\lambda x^\rho )^2=x^\lambda\cdot xx^\lambda,\qquad
(x^{\lambda^2}\cdot x^\lambda x^\rho )^\lambda =(x^{\lambda^2}\cdot
x^\lambda x^\lambda )^\rho,
\end{displaymath}
\begin{displaymath}
(x^\lambda\cdot xx)^{\lambda^2}\cdot (x^\lambda\cdot xx)^\lambda
(x^\lambda\cdot xx)^\rho =x^\lambda\cdot xx^\lambda.
\end{displaymath}
\item
\begin{displaymath}
(x\cdot x^\rho y^\rho )^\lambda=(x\cdot x^\rho y^\lambda
)^\rho,\qquad (x^\lambda\cdot xy^{\rho^2})^\lambda =(x^\lambda\cdot
xy)^\rho,
\end{displaymath}
\begin{displaymath}
(x^{\lambda^2}\cdot x^\lambda y^\rho )^\lambda =(x^{\lambda^2}\cdot
x^\lambda y^\lambda )^\rho,\qquad (x^\lambda\cdot xy)^\lambda
=(x^\lambda\cdot xy^{\lambda^2})^\rho.
\end{displaymath}
\item $|J_\lambda |=2$ iff $|J_\rho |=2$ iff $J_\lambda =J_\rho$ iff
$L$ is a LSIPL iff RSIPL
\end{enumerate}
\end{mylem}
{\bf Proof}\\
The whole these is gotten by intuitive use of (b), (c) and (d) of
Theorem~\ref{3:16}.

\begin{mycor}\label{3:17.1}
Let $L$ be a CC-loop. The following are equivalent.
\begin{enumerate}
\item $L$ is a power associativity loop
\item $L$ is a 3-PAPL.
\item $L$ obeys $x^\rho =x^\lambda$ for all $x\in L$.
\item $L$ is a LSIPL.
\item $L$ is a RSIPL.
\end{enumerate}
\end{mycor}
{\bf Proof}\\
The proof the equivalence of the first three is shown in Lemma~3.20
of \cite{phd78} and mentioned in Lemma~1.2 of \cite{phd151}. The
proof of the equivalence of the last two and the first three can be
deduced from the last result of Lemma~\ref{3:17}.

\begin{myrem}
This new algebraic definition gives more insight into the algebraic
properties of Osborn loop. Particularly, it can be used to fine tune
some recent equations on CC-loop as shown in works of Kunen, Kinyon,
Phillips and Drapal; \cite{phd35,phd36,phd47}, \cite{phd37,phd38},
 \cite{phd78}. In fact, in \cite{phd35,phd78}, the authors focussed on the mapping
$E_x=R_xR_{x^\rho}=\theta_x^{-1}$ where $\theta_x =L_xL_{x^\lambda}$
and were able to established study its algebraic properties in a
CC-loop. So we can see that the investigations of $E_x$ in CC-loops
by Kunen, Kinyon and Phillips is a bit in line with what Huthnance
\cite{phd44} did with $\theta_x$ in a universal WIPL and WIP Osborn
loop. In this work, attention has been paid primarily on Osborn
loops. So this study is a general overview of the earlier ones. The
identities LSIPL and RSIPL are appearing for the first time.
\end{myrem}

\subsection{Application Of An Osborn Loop Identity To Cryptography}
Among the few identities that have been established for Osborn loops
in Theorem~\ref{3:16}, we would recommend one of them for
cryptography in a similar spirit in which the cross inverse property
has been used by Keedwell \cite{phd176}. It will be recalled that
CIPLs have been found appropriate for cryptography because of the
fact that the left and right inverses $x^\lambda$ and $x^\rho$ of an
element $x$ do not coincide unlike in left and right inverse
property loops, hence this gave rise to what is called 'cycle of
inverses' or 'inverse cycles' or simply 'cycles' i.e finite sequence
of elements $x_1,x_2,\cdots ,x_n$ such that $x_k^\rho
=x_{k+1}~\bmod{n}$. The number $n$ is called the length of the
cycle. The origin of the idea of cycles can be traced back to Artzy
\cite{phd140,phd193} where he also found there existence in WIPLs
apart form CIPLs. In his two papers, he proved some results on
possibilities for the values of $n$ and for the number $m$ of cycles
of length $n$ for WIPLs and especially CIPLs. We call these "Cycle
Theorems" for now.

\paragraph{}
In the course of this study(Lemma~\ref{3:17}), it has been
established that in an Osborn loop, $J_\lambda =J_\rho$, LSIP and
RSIP are equivalent conditions. Therefore, in a CC-loop, the power
associativity property, 3-PAPL, $x^\rho =x^\lambda$, LSIP and RSIPL
are equivalent. Thus, Osborn loops without the LSIP or RSIP will
have cycles(even long ones). This exempts groups, extra loops, and
Moufang loops but includes CC-loops, VD-loops and universal WIPLs.
Precisely speaking, non-power associative CC-loops will have cycles.
So broadly speaking and following some of the identities in
Lemma~\ref{3:17}, Osborn loops that do not have the LSIP or RSIP or
3-PAPL or weaker forms of inverse property, power associativity and
diassociativity to mention a few, will have cycles(even long ones).
The next step now is to be able to identify suitably chosen
identities in Osborn loops, that will do the job the identity
$xy\cdot x^\rho =y$ or its equivalents does in the application of
CIPQ to cryptography. These identities will be called Osborn
cryptographic identities(or just cryptographic identities).
\begin{mydef}\label{5:1}
Let $\mathcal{Q}=(Q, \cdot ,\backslash ,/)$ be a quasigroup. An
identity $w_1(x,x_1,x_2,x_3,\cdots )=w_2(x,x_1,x_2,x_3,\cdots )$
where $x\in Q$ is fixed, $x_1,x_2,x_3,\cdots\in Q$, $x\not\in
\{x_1,x_2,x_3,\cdots\}$ is said to be a cryptographic identity(CI)
of the loop $\mathcal{Q}$ if it can be written in a functional form
$xF(x_1,x_2,x_3,\cdots )=x$ such that $F(x_1,x_2,x_3,\cdots )\in
\mathcal{M}\textrm{ult}(\mathcal{Q})$. $F(x_1,x_2,x_3,\cdots )=F_x$
is called the corresponding cryptographic functional(CF) of the CI
at $x$.
\end{mydef}
\begin{mylem}\label{5:2}
Let $\mathcal{Q}=(Q, \cdot ,\backslash ,/)$ be a loop with identity
element $e$ and let $CF_x(\mathcal{Q})$ be the set of all CFs in
$\mathcal{Q}$ at $x\in Q$. Then, $CF_x(\mathcal{Q})\le
\mathcal{M}\textrm{ult}(\mathcal{Q})$ and $CF_e(\mathcal{Q})\le
\textrm{Inn}(\mathcal{Q})$.
\end{mylem}
\begin{mylem}\label{5:3}
Let $\mathcal{Q}=(Q, \cdot ,\backslash ,/)$ be a quasigroup.
\begin{enumerate}
\item $T_{(x)}\in CF_y(\mathcal{Q})$ if and only if $y\in C(x)$,
\item $R_{(x,y)}\in CF_z(\mathcal{Q})$ if and only if $z\in N_\lambda (x,y)$,
\item $L_{(x,y)}\in CF_z(\mathcal{Q})$ if and only if $z\in N_\rho (x,y)$,
\end{enumerate}
where $N_\lambda (x,y)=\{z\in Q|zx\cdot y=z\cdot xy\}$, $N_\rho
(x,y)=\{z\in Q|y\cdot xz=yx\cdot z\}$ and $C(x)=\{y\in Q|xy=yx\}$.
\end{mylem}
\begin{mylem}\label{5:4}
Let $\mathcal{Q}=(Q, \cdot ,\backslash ,/)$ be an Osborn loop with
identity element $e$. Then, the identity $yx=x(x^\lambda y\cdot x)$
is a CI with its CF $F_e\in CF_e(\mathcal{Q})$.
\end{mylem}
\begin{myrem}
The identity $yx=x(x^\lambda y\cdot x)\Leftrightarrow y=[x(x^\lambda
y\cdot x)]/x$ is more "advanced" than the CIPI and hence will posse
more challenge for an attacker(even than the CIPI) to break into a
systems. As described by Keedwell, for a CIP, it is assumed that the
message to be transmitted can be represented as single element $x$
of a CIP quasigroup and that this is enciphered by multiplying by
another element $y$ of the CIPQ so that the encoded message is $yx$.
At the receiving end, the message is deciphered by multiplying by
the inverse of $y$. But for the identity $y=[x(x^\lambda y\cdot
x)]/x$, procedures of enciphering and deciphering are more than one
in an Osborn loop.
\end{myrem}

\end{document}